\documentclass{amsart}
\usepackage{hyperref}
\usepackage[english]{babel}
\usepackage{enumerate}
\usepackage{mathrsfs}
\usepackage{amsmath}
\usepackage{amssymb}
\usepackage{amsthm}
\newcommand{\al}[1]{\mathcal{#1}}

\newcommand{\ind}{\mbox{ind}}

\newcommand{\mb}[1]{\mathbb{#1}}
\newcommand{\Img}{\mbox{Im}}

\newcommand{\ran}{\mbox{ran}}
\newcommand{\R}{\mathbb{R}}

\newcommand{\C}{\mb{C}}
\theoremstyle{definition}

\newtheorem{remark}{Remark}[section]

\newtheorem{lemma}{Lemma}[section]
\newtheorem{proposition}{Proposition}[section]
\newtheorem{example}{Example}[section]

\numberwithin{equation}{section}
\begin{document}
\title[On the connected components of the conjugacy class of
projectors]{On the connected components of the conjugacy class %
of projectors on $ \ell_p\oplus\ell_q $}
\thanks{This work was supported by INHA %
UNIVERSITY Research Grant}
\author{Daniele Garrisi}
\email{daniele.garrisi@inha.ac.kr}
\address{West Building, Office No. 5W443\\
Department of Mathematics Education\\
Inha University\\
253 Yonghyun-Dong, Nam-Gu\\
Incheon, South Korea 402-751}
\begin{abstract}
We characterize the projectors $ P $ on a Banach space $ E $ 
having the property of being connected to all the others projectors 
obtained as a conjugation of $ P $. 
Using this characterization we show an example of Banach 
space where the conjugacy class of a projector splits into several 
path-connected components, and describe the conjugacy classes
of projectors onto subspaces of $ \ell_p\oplus\ell_q $ with
$ p\neq q $.
\end{abstract}
\maketitle
\section{Introduction}
Let $ \al{A} $ be a Banach algebra and denote by 
$ \mathcal{P}(\al{A}) $ the space of \textsl{projectors}, that is
the set of the elements $ p\in\al{A} $ 
such that $ p^2 = p $, endowed with the subspace topology.
We denote by $ G(\al{A}) $ the group of invertible elements of 
$ \al{A} $. Given a topological space $ X $, we use the notation 
$ \pi_0 (X) $ for set of path-connected components of $ X $, and 
$ [x] $ for the component which contains $ x $.

Given two projectors $ p $ and $ q $ it is possible to define
another relation
\begin{equation*}
pRq:\exists g\in G(\al{A}) \text{ such that } gp = qg,
\quad R_p := \{q\mid pRq\}.
\end{equation*}
If $ pRq $, we also say that $ p $ is conjugated to $ q $, and
$ R_p $ is the conjugacy class of $ p $.
In general, for every Banach algebra, the inclusion 
$ [p]\subseteq R_p $ holds. A proof of this fact can be found in 
\cite[Proposition~4.2]{PR87}, or 
\cite[Ch. I.4 \S 6]{Kat95}.
Actually, the following fact holds true: if $ [p] = [q] $, 
there exists a continuous path $ u\colon [0,1]\to G(\al{A}) $ such that 
\begin{equation}
\label{eq.pairs-of-projectors}
u(1) p = q u(1)\text{ and } u(0) = 1.
\end{equation}
In \cite[7.13.1]{PR87}, G.~Porta and L.~Recht provided an example of 
Banach algebra and a projector $ p $ such that $ [p]\subsetneq R_p $.
The algebra $ \al{A} $ is $ C(S^3,M(2,\C)) $ and the 
projector is defined as follows: if $ x := (x_1,x_2,x_3,x_4) $ is 
a point of $ S^3 $, then
\begin{equation}
\label{eq.pr}
p(x) := 
\begin{pmatrix}
  \overline{z}_0 z_0 &  \overline{z}_1 z_0 \\
  \overline{z}_0 z_1 &   \overline{z}_1 z_1
\end{pmatrix},\quad q(x) := 1 - p(x)
\end{equation}
where $ z_0 := x_1 + ix_2 $ and $ z_1 := x_3 + ix_4 $. 
Therefore, we wonder whether a similar example may
occur in algebras of bounded operators on a Banach space $ E $, 
for a suitable choice of $ E $. That is, whether there are two
projectors $ P $ and $ Q $ in $ \al{P}(\al{L}(E)) $ which are conjugated
to each other, but $ [P]\neq [Q] $, for a suitable choice of $ E $.
This happens in the space $ E := \ell_p\oplus\ell_q $ with $ 1 \leq p,q $
and $ p\neq q $.

In \S\ref{sect.necessary-sufficient}, we provide necessary and 
sufficient conditions to a projector 
$ P\in\mathcal{P}(\al{L}(E)) $ ensuring
the equality $ R_P = [P] $. This conditions is expressed in terms
of the connected components of the general linear group of
$ \ker(P),\ran(P) $ and $ E $. In 
\S\ref{sect.explicit-construction}, 
we address projectors $ P $ such that $ \ran(P) $ can
be obtained as direct sum of a subspace of $ \ell_p\oplus\{0\} $
and a subspace of $ \{0\}\oplus\ell_q $. 
We show that there are two possible behaviours: either $ R_P $ is 
path-connected or $ R_P $ has infinitely many connected components. 
\section{A criterion to establish whether 
$ R_P = [P] $}
\label{sect.necessary-sufficient}
Let $ P $ be a projector of $ \al{L}(E) $. We set $ X := \ran(P) $ 
and $ Y := \ker(P) $.
\begin{proposition}
\label{prop.conjugation}
Given a Banach space $ E $, $ R_P = [P] $ if and only if 
\begin{equation*}
\varphi_*\colon\pi_0 (GL(X)\times GL(Y))\to \pi_0 (GL(E)),\quad
\varphi(A,B) := A\oplus B
\end{equation*}
is surjective.
\end{proposition}
\begin{proof}
Suppose that $ R_P \subseteq [P] $. Given $ T\in GL(E) $, we 
set
\begin{equation}
\label{eq.2}
Q := TPT^{-1}.
\end{equation}
Since $ Q $ is a projector, it is path-connected to $ P $.
Then, from \eqref{eq.pairs-of-projectors}, there exists
a path of invertible operators $ U $ such that
\begin{equation}
\label{eq.3}
SP = QS,\quad S := U(1).
\end{equation}
From \eqref{eq.2} and \eqref{eq.3}, it follows that
$ S^{-1} T $ commutes with $ P $ or, equivalently, $ S^{-1} T $ 
is in the image of $ \varphi $. Hence $ S^{-1} T\in\Img(\varphi) $
while, from the definition of $ U(t) $, we have $ [S^{-1} T] = [T] $.
Then, $ [T]\in\Img(\varphi_*) $. Conversely, suppose that 
$ \varphi_* $ is surjective and consider $ TPT^{-1} $, a conjugation 
of $ P $. Let $ (A,B) $ be a pair such that 
$ [A\oplus B] = [T^{-1}] $. Since $ A\oplus B $ commutes with $ P $
we have
\begin{equation*}
(T(A\oplus B)) P (T (A\oplus B))^{-1} = TPT^{-1},
\end{equation*}
while $ [T(A\oplus B)] = [id_E] $. Then $ [TPT^{-1}] = [P] $.
\end{proof}
\begin{remark}
\label{rem.1}
Actually, we proved that $ TPT^{-1} $ belongs to the same
connected component of $ P $ if and only if $ [T] $ belongs
to the image of $ \varphi_* $.
\end{remark}
\begin{example}
When $ GL(E) $ is path-connected, $ R_P = [P] $. In fact,
given $ U\in GL(E) $ and a path $ U(t) $ which connects
$ U $ to $ id_E $, the path $ U(t) P U(t)^{-1} $ connects
$ P $ to $ UPU^{-1} $.
\end{example}
\begin{example}
Considering cases where the linear group is not
connected does not automatically imply the existence
of a projector where $ [P]\subsetneq R_P $. For instance, if
$ E = \R^n $, the linear group is $ GL(n,\R) $ which
consists of two connected components characterized by the
sign of the determinant; given $ T\in GL(n,\R) $ such that
$ \det(T) < 0 $ it is possible to choose $ A\in GL(X) $
and $ B\in GL(Y) $ such that $ \det(A)\cdot\det(B) < 0 $.
Therefore $ \varphi_* ([A\oplus B]) = [T] $, and $ \varphi_* $
is surjective.
\end{example}
\section{The conjugacy class in 
$ \ell_p\oplus\ell_q $}
\label{sect.explicit-construction}
We will use a construction of A.~Douady, \cite{Dou65}, devised
with the purpose of showing the existence of infinite-dimensional
Banach spaces whose linear group is not path-connected.

Let $ E = F\oplus G $ be Banach spaces such that $ F $ and 
$ G $ are isomorphic to their closed subspaces of co-dimension one
and such that $ \al{L}(F,G) = \al{L}_{SS} (F,G) $ (the subspace of
strictly singular operators).
We set 
\begin{equation*}
P_F (x,y) := (x,0),\quad P_G (x,y) = (0,y).
\end{equation*}
Given an operator $ A\in\al{L}(F) $, we use the notation
$ \ind(A) $ for the Fredholm index, whenever appropriate.
\begin{lemma}[From {\cite[Proposition~1]{Dou65}}]
\label{lem.douady}
There exists a continuous, surjective group homomorphism 
$ j\colon GL(E)\rightarrow\mb{Z} $
such that $ j(T) = \ind(P_F T P_F) $. Moreover,
\begin{enumerate}[\ \ \ (a)]
\item\label{item.a}
 $ \ind(P_F T P_F) + \ind(P_G T P_G) = 0 $
\item if $ GL(F) $ and $ GL(G) $ are path-connected, 
$ j $ is injective. 
\end{enumerate}
\end{lemma}
For the proof of the fact that $ j $ is a group homomorphism and
that is surjective, we refer to Lemma~1 and Lemma~2 of \cite{Mit70}.
An explicit invertible operator $ V $ such that $ j(V) = 1 $ is the
example provided in \cite[p.~6]{Mit70}. Let $ F_1 $ and $ G_1 $ be
two closed subspaces of co-dimension one in $ F $ and $ G $, respectively
and $ \sigma\colon F_1\to F $ and $ \tau\colon G_1\to G $ be
two isomorphisms. Let also $ v $ and $ w $ be such that
\begin{equation*}
F = F_1\oplus\langle v\rangle,\quad G = G_1\oplus\langle w\rangle.
\end{equation*}
We define
\begin{equation}
\label{eq.j=1}
V(x + s v,y) = (\sigma(x),s w + \tau^{-1} (y)).
\end{equation}
where $ \ind(P_F V P_F) = 1 $. For the injectivity of the map $ j $
when $ GL(F) $ and $ GL(G) $ are path connected, we
refer to \cite{Dou65} and \cite{Jan65}.

Hereafter, we set $ F := \ell_p $ and $ G := \ell_q $.
Let $ P $ be such that $ \ran(P) = X_1\oplus X_2 $, where
$ X_1\subseteq\ell_p $ and $ X_2\subseteq\ell_q $ are two closed 
and complemented subspaces. Let $ Y_1 $ and $ Y_2 $
be closed subspaces such that
\begin{equation*}
X_1 \oplus Y_1 = \ell_p,\quad X_2\oplus Y_2 = \ell_q.
\end{equation*}
\begin{proposition}
If both $ Y_1 $ and $ Y_2 $, or both $ X_1 $ and $ X_2 $ have infinite
dimension, then $ [P] = R_P $. 
That is, any projector conjugated to $ P $ is also path-connected to it.
Otherwise, $ \pi_0 (R_P)\approx\mb{Z} $, that is in the conjugacy class
of $ P $ there are infinitely many connected components.
\end{proposition}
\begin{proof}
From \cite[Theorem~2.24, p. 35]{AK06}, closed infinite-dimensional 
and complemented subspaces of $ \ell_p $ and $ \ell_q $ are isomorphic to 
$ \ell_p $ and $ \ell_q $, respectively. 
Moreover, from \cite{Neu67}, the linear groups of
$ \ell_p $ and $ \ell_q $ are contractible, and 
$ \al{L}_{SS} (\ell_p,\ell_q) = \al{L} (\ell_p,\ell_q) $, from
\cite[Pitt's~Theorem, p. 32]{AK06} and 
\cite[Theorem~2.1.9, p. 33]{AK06}. Therefore, the decomposition 
\begin{equation*}
(X_1\oplus Y_1) \bigoplus (X_2 \oplus Y_2) = E,\quad
F := X_1\oplus Y_1,\quad G := X_2\oplus Y_2
\end{equation*}
satisfies all the assumptions of Lemma~\ref{lem.douady}, making $ j $ a 
group isomorphism on $ \pi_0 (GL(E)) $. From 
Proposition~\ref{prop.conjugation}, we need to check whether the inclusion 
\begin{equation*}
\varphi_*\colon
\pi_0 (GL(X_1\oplus X_2)\times GL(Y_1\oplus Y_2))\to \pi_0 (GL(E))
\end{equation*}
is surjective. In particular, we want to know
whether we can obtain an invertible operator in $ T $ on $ E $ such that 
$ j(T) = 1 $. Let
\begin{equation*}
(A,B)\in GL(X_1\oplus X_2)\times GL(Y_1\oplus Y_2)
\end{equation*}
be an arbitrary element and $ T := A\oplus B $.
We denote by $ P_1 $ and $ P_2 $ the projectors on
$ X_1 $ and $ X_2 $, and by $ Q_1 $ and $ Q_2 $, the projectors
on $ Y_1 $ and $ Y_2 $. We have
\begin{align}
\label{eq.system-1}
1 &= \ind(P_1 AP_1) + \ind(Q_1 B Q_1)\\
\label{eq.system-2} 
0 &= \ind(Q_1 B Q_1) + \ind(Q_2 B Q_2)\\
\label{eq.system-3} 
0 &= \ind(P_1 A P_1) + \ind(P_2 A P_2).
\end{align}
The first equality comes from the requirement $ j(T) = 1 $.
The second and the third equalities follows from 
\eqref{item.a} of Lemma~\ref{lem.douady} when it is applied to
the spaces $ Y_1,Y_2 $ and $ X_1,X_2 $.

In finite-dimensional spaces, all the operators are
Fredholm with index zero. Suppose that the assumption of the proposition
is not fulfilled, that is, for instance, $ X_1 $ and $ Y_1 $ have finite
dimension. Therefore, in \eqref{eq.system-1}, we obtain $ 1 = 0 + 0 $;
if $ X_1 $ and $ Y_2 $ have finite dimension, therefore, from 
\eqref{eq.system-2}, we obtain
\begin{equation*}
\ind(Q_2 B Q_2) = 0 = \ind(Q_1 B Q_1)
\end{equation*}
which yields to
\begin{equation*}
1 = \ind(P_1 AP_1) = 0
\end{equation*}
combining \eqref{eq.system-1} and the fact that $ X_1 $ has finite 
dimension. A similar arguments gives a contradiction when $ X_2 $ and 
$ Y_1 $ both have finite dimension.

Now, we suppose that both $ X_1 $ and $ X_2 $ have infinite dimension.
Then, we can define on $ X_1\oplus X_2 $ an invertible operator $ V $ 
similar to \eqref{eq.j=1}. Then
\begin{equation*}
j(\varphi(V,id_{Y_1\oplus Y_2})) = 1.
\end{equation*}
And, if both $ Y_1 $ and $ Y_2 $ have infinite dimension,
\begin{equation*}
j(\varphi(id_{X_1\oplus X_2},V)) = 1.
\end{equation*}
\end{proof}
Therefore, in $ \ell_p\oplus\ell_q $, the projector
$ P_F $ (or $ P_G $) is an example
of projector such that $ [P]\neq R_P $. We notice that $ P_F $
is not conjugated to $ id_E - P_F = P_G $ ($ \ell_p $ is
not isomorphic to $ \ell_q $). This might suggest that the 
features of 
the example in \eqref{eq.pr} were not entirely transposed to the 
setting of linear operators. The following remark explains why 
this is not possible
\begin{remark}
\label{rem.1-p}
If $ P $ is in $ \al{P}(\al{L}(E)) $, and $ id_E - P\in R_P $,
then $ [id_E - P] = [P] $.
\end{remark}
We conclude this paper by showing that our example is substantially 
different from the one obtained in \cite{PR87}. We mean that 
$ C(S^3,M(2,\C)) $ does not have the algebraic structure of the space of 
bounded operators on a Banach space $ E $. This would already follow
from Remark~\ref{rem.1-p}, but we prefer to give a proof which is 
independent from the content of the work in \cite{PR87}.
\begin{proposition}
\label{prop.ckM}
Given a non-empty compact space $ K $, a natural number $ n $ and a Banach 
space $ E $, $ C(K,M(n,\C)) $ is not isomorphic to $ \al{L}(E) $, unless 
$ E\simeq\C^n $ and $ K $ is a singleton.
\end{proposition}
\begin{proof}
We look at the centers of the two algebras. We have
\begin{equation*}
Z(\al{L}(E)) = \langle id_E \rangle,\quad
Z(C(K,M(n,\C))) = C(K,\langle id_n \rangle).
\end{equation*}
If $ C(K,M(n,\C))\simeq\al{L}(E) $ there is an algebra isomorphism between
the two centers. Therefore, $ Z(\al{L}(E))\simeq Z(C(K,M(n,\C))) $.
The first space has dimension one, while the second space has infinite 
dimension unless $ K $ is finite. If $ K $ is a singleton, then the 
dimension is one. In that case $ C(K,M(n,\C))\simeq M(n,\C) $. In 
conclusion, $ M(n,\C) $ is isomorphic to $ \al{L}(E) $ only in one case: 
when $ E \simeq\C^n $.
\end{proof}

\providecommand{\bysame}{\leavevmode\hbox to3em{\hrulefill}\thinspace}
\providecommand{\MR}{\relax\ifhmode\unskip\space\fi MR }
\providecommand{\MRhref}[2]{%
  \href{http://www.ams.org/mathscinet-getitem?mr=#1}{#2}
}
\providecommand{\href}[2]{#2}

\end{document}